# ON DISCRIMINATING BETWEEN LONG-RANGE DEPENDENCE AND CHANGES IN MEAN[1]

By István Berkes,[2] Lajos Horváth, Piotr Kokoszka[3] and Qi-Man Shao[4]

*Graz University of Technology and A. Rényi Institute of Mathematics, University of Utah, Utah State University, and University of Oregon and Hong Kong University of Science and Technology*

We develop a testing procedure for distinguishing between a long-range dependent time series and a weakly dependent time series with change-points in the mean. In the simplest case, under the null hypothesis the time series is weakly dependent with one change in mean at an unknown point, and under the alternative it is long-range dependent. We compute the CUSUM statistic $T_n$, which allows us to construct an estimator $\hat{k}$ of a change-point. We then compute the statistic $T_{n,1}$ based on the observations up to time $\hat{k}$ and the statistic $T_{n,2}$ based on the observations after time $\hat{k}$. The statistic $M_n = \max[T_{n,1}, T_{n,2}]$ converges to a well-known distribution under the null, but diverges to infinity if the observations exhibit long-range dependence. The theory is illustrated by examples and an application to the returns of the Dow Jones index.

## 1. Introduction.

The present paper develops a testing procedure for distinguishing between a long-range dependent time series and a weakly dependent time series with change-points in the mean.

Many geophysical time series records have long been known to exhibit long nonperiodic cycles or persistent deviations from the mean. In the mid-1960s Mandelbrot and his collaborators proposed the use of self-similar processes, most notably fractional Brownian motion, to model such records;

Received July 2003; revised July 2005.
[1]Supported by NSF Grant INT-0223262 and NATO Grant PST.EAP.CLG 980599.
[2]Supported in part by OTKA Grants T 037886 and T 043037.
[3]Supported in part by NSF Grant DMS-04-13653.
[4]Supported in part by NSF Grant DMS-01-03487, and Grants R-146-000-038-101 and R-1555-000-035-112 at the National University of Singapore.

see, for example, [33]. Over a decade later, Granger and Joyeux [21] and Hosking [25] (see also [1]) introduced fractional ARIMA processes, which are approximately self-similar and offer a much greater modeling flexibility. If their fractional differencing parameter $d$ satisfies $0 < d < 1/2$, these processes are stationary and possess long-range dependence, or long memory, in the sense that the autocovariance function is not absolutely summable (decays like $k^{2d-1}$, as the lag $k \to \infty$). In the 1980s there was substantial interest in using long memory processes to model macroeconomic time series, whereas, in the 1990s, the focus shifted to modeling the volatility of returns on speculative assets by such processes; an in depth discussion and relevant references are provided in [23]. Following the pioneering work of Leland et al. [30] and Paxson and Floyd [38], self-similar processes have also increasingly been used to model certain aspects of computer network traffic; see [36]. There are many other fields where models exhibiting long-range dependence have been used; see [14] for a recent extensive review.

Even though modeling certain time series in the aforementioned fields by means of long-range dependent processes has become quite widespread, especially in geophysics, it is clear that a series with long periods where the observations are away from the mean can also naturally be modeled by a nonstationary process whose mean changes. Bhattacharya, Gupta and Waymire [9] used mathematical arguments to show that the so-called *Hurst effect*, which motivated Mandelbrot and his collaborators to advocate the use of self-similar processes, can also be explained if the observations $X_k$ are assumed to follow the model $X_k = Y_k + f(k)$, where $Y_k$ is a weakly dependent stationary process and $f$ is a deterministic function. That research was elaborated on by Giraitis, Kokoszka and Leipus [16] who showed that several statistics akin to the modified $R/S$ statistic of Lo [31] diverge to infinity under either long-range dependence or weak dependence with change-points. In a similar spirit, Diebold and Inoue [13] argued that the appearance of long memory can be explained by some econometric models which involve changes in their defining parameters. Mikosch and Stărică [34, 35] asserted that what had been seen by many as long memory in the volatility of returns is, in fact, a manifestation of changes in the parameters of the underlying GARCH-type models. In the context of network traffic, similar findings are reported in [26]. The above list of references is not exhaustive, but it emphasizes that it is difficult to distinguish a truly long-range dependent process from a process with some form of nonstationarity, including shifts in mean. Standard tools like ACF plots and periodogram-based spectral estimates behave in a very similar way under these two alternatives. There are also a number of long memory tests designed to test the null hypothesis of weak dependence against an alternative of long-range dependence and change-point tests developed to test the same null hypothesis but against a change-point alternative. Most long memory tests reject in the presence of change-points and many change-point tests reject in the presence of long memory.

The answer to the question of which approach to use will often depend on a specific application at hand. A long-range dependent process may, for example, provide a parsimonious description of a long, possibly nonstationary, time series. On the other hand, to construct short term forecasts of a possibly self-similar process, it might be advisable to fit an ARMA model to the most recent stretch of data after the last estimated change-point. In many applications, however, such as, for example, constructing long term forecasts, it does matter which model better fits the data. We refer to [10] and [28] for some relevant financial applications. Formal statistical tests which would help decide if a particular time series is better described as a realization of a long-range dependent process or as a realization of a weakly dependent process with change-points are therefore of value. There has, however, not been much research in this direction. Künsch [29] proposed a periodogram based procedure to discriminate between a long-range dependent process and the process $X_k = Y_k + f(k)$ with a monotonic function $f$ and Gaussian weakly dependent $Y_k$. Heyde and Dai [24] showed that procedures for detecting long memory which are based on a smoothed periodogram are robust in the presence of small trends. These ideas were recently developed by Sibbertsen and Venetis [42] who proposed a test based on a difference between the Geweke and Porter-Hudak [15] estimator of $d$ and its version based on the tapered periodogram.

A main objective of the present paper is to develop the theory underlying a test procedure for discriminating between long-range dependence and weak dependence with change-points in mean. The proposed test is a simple time domain procedure based on a CUSUM statistic for the partial sums, which is perhaps the most extensively used statistic for detecting and estimating change-points in mean. To describe the idea, suppose that, under the null hypothesis, the time series is weakly dependent with one change in mean and under the alternative, it is long-range dependent. Consider the CUSUM statistic $T_n$ defined by (3.1). Using $T_n$, we can construct an estimator $\hat{k}$ of the change-point (no matter if a change-point exists or not). We then compute the statistic $T_{n,1}$ based on the observations up to time $\hat{k}$ and the statistic $T_{n,2}$ based on the observations after time $\hat{k}$. The statistic $M_n = \max[T_{n,1}, T_{n,2}]$ converges to a well-known distribution under the null (cf. Corollary 2.1), but diverges to infinity under the alternative.

Our theory uses the almost sure asymptotics for the Bartlett variance estimator $s_n^2$ stated in Theorem A.1 which was established in [8]. For a weakly dependent process, $s_n^2$ is an estimator of the variance of the sample mean or of the spectral density at frequency zero. Estimators of this type have been extensively studied in the time series literature in the last half century and go back to the work of Bartlett [5], Grenander and Rosenblatt [22] and Parzen [37]. Andrews [3] provides a more recent perspective. As far as we know, all consistency results pertaining to the class of kernel estimators such as $s_n^2$ establish convergence in an $L^p$ norm or in probability. Such results might possibly be applied in our context after some additional technical work, but we are not aware of any convergence in probability results which would allow us to establish our main results, Theorems 2.1 and 2.2, under weaker conditions. Moreover, almost sure convergence offers a convenient approach based on the observation that if $Z_n \stackrel{\text{a.s.}}{\to} 0$ and $k_n \stackrel{P}{\to} \infty$, then $Z_{k_n} \stackrel{P}{\to} 0$ [see, e.g., the argument justifying (B.11)].

The paper is organized as follows. In Section 2 we formulate the assumptions, describe the testing procedure in a simple illustrative situation and state the relevant theorems. Section 3 discusses the broader applicability of the procedure, provides some additional background and examples and concludes with an application to returns of the Dow Jones index. The appendices contain the proofs.

**2. Assumptions and the testing procedure.** To focus attention and lighten the notation, we concentrate in this section on a situation where the observations can either follow a model with one change in the mean of weakly dependent time series or are long-range dependent. In Section 3 we explain how the proposed procedure can be used in a situation when there is an upper bound on the number of possible changes in the mean.

The observations $X_i$ follow a change-point model if

$$X_i = \begin{cases} \mu + Y_i, & 1 \le i \le k^*, \\ \mu + \Delta + Y_i, & k^* < i \le n. \end{cases} \tag{2.1}$$

In (2.1) $k^*$ is the unknown time of a possible change in mean, and the means $\mu$ and $\mu + \Delta$ are also unknown. The sequence $\{Y_i\}$ is assumed to have mean zero and to be weakly stationary in a sense made precise by Assumption 2.1. Recall that, for a fourth-order stationary sequence $\{Y_k\}$ with mean 0 and $\gamma_j = \mathrm{Cov}(Y_0, Y_j)$, the fourth-order cumulant is defined by

$$\kappa(h,r,s) = E[Y_k Y_{k+h} Y_{k+r} Y_{k+s}] - (\gamma_h \gamma_{r-s} + \gamma_r \gamma_{h-s} + \gamma_s \gamma_{h-r}). \tag{2.2}$$

Assumption 2.1. The sequence $\{Y_k\}$ is fourth-order stationary with mean 0 and autocovariance function $\gamma_j = \mathrm{Cov}(Y_0, Y_j)$, and the following conditions hold:

$$n^{-1/2} \sum_{1 \le j \le nt} Y_j \stackrel{d}{\to} \sigma W(t) \qquad \text{in } D[0,1] \tag{2.3}$$

for some $\sigma > 0$ and

$$\sum_j |\gamma_j| < \infty, \tag{2.4}$$

$$\sup_h \sum_{r,s} |\kappa(h,r,s)| < \infty. \tag{2.5}$$

REMARK 2.1. By the Skorokhod–Wichura–Dudley representation (see, e.g., [41]), condition (2.3) is equivalent to the following condition: There are Wiener processes $W_n(t), t \in [0,1]$, such that

$$\sup_{0\le t\le 1}\left| n^{-1/2} \sum_{1\le j\le nt} Y_j - \sigma W_n(t)\right| = o_P(1). \tag{2.6}$$

Condition (2.6) is often more convenient to refer to in the proofs.

We now make precise the statement that the observations $\{X_i\}$ are long-range dependent. In the following $W_H(t)$ stands for the fractional Brownian motion with parameter $H$, that is, a Gaussian process with mean zero and covariances

$$E[W_H(t)W_H(s)] = (t^{2H} + s^{2H} - |t-s|^{2H})/2.$$

If $1/2 < H < 1$, the increments of the fractional Brownian motion are long-range dependent. It is convenient to identify the self-similarity parameter $H$ with the differencing parameter $d$ introduced in Section 1 via the relation $H = d + 1/2$ because the increments of $W_H$, which form a stationary process, have the same rate of decay of the autocovariance function as a fractional ARIMA with $d = H - 1/2$; see, for example, Section 7.13 of [40]. In condition (2.9) of Assumption 2.2 below, and throughout the paper, $a_j \sim b_j$ means that $\lim_{j\to\infty} a_j/b_j = 1$.

ASSUMPTION 2.2. The sequence $\{X_j\}$ is fourth-order stationary with $\mu = EX_j$ and $\gamma_j = \mathrm{Cov}(X_0, X_j)$ and satisfies the following conditions:

$$\frac{1}{n^H} \sum_{1\le j\le nt} (X_j - \mu) \xrightarrow{d} c_H W_H(t) \qquad \text{in } D[0,1] \tag{2.7}$$

for some $c_H > 0$ and

$$\tfrac{1}{2} < H < 1. \tag{2.8}$$

Moreover,

$$\gamma_j \sim c_0 j^{2H-2} \tag{2.9}$$

for some $c_0 > 0$, and the cumulants (2.2) satisfy

$$\sup_h \sum_{-n\le r,s\le n} |\kappa(h,r,s)| = O(n^{2H-1}). \tag{2.10}$$

The cumulant condition (2.5) is weaker than the traditional condition in which $\sup_h$ is replaced by $\sum_h$; see, for example, [2] and [3]. Condition (2.10) is a natural counterpart of (2.5) and holds for the extensively used fractional ARIMA models. For these models, the range (2.8) corresponds to $0 < d < 1/2$. We do not consider $-1/2 < d < 0$ because realizations of such processes do not exhibit apparent shifts in mean.

We wish to test

H$_0$: The observations $X_1, \dots, X_n$ follow the change point model (2.1) with the $Y_i$ satisfying Assumption 2.1

against

H$_A$: The observations $X_1, \dots, X_n$ are long-range dependent, that is, satisfy Assumption 2.2.

In order to define the test statistic, we first introduce a change-point estimator,

$$\hat{k} = \min\left\{k : \max_{1\le i\le n}\left|\sum_{1\le j\le i} X_j - \frac{i}{n}\sum_{1\le j\le n} X_j\right| = \left|\sum_{1\le j\le k} X_j - \frac{k}{n}\sum_{1\le j\le n} X_j\right|\right\}. \tag{2.11}$$

Next we define the statistics

$$T_{n,1} = \frac{1}{s_{n,1}}\hat{k}^{-1/2}\max_{1\le k\le \hat{k}}\left|\sum_{1\le i\le k} X_i - \frac{k}{\hat{k}}\sum_{1\le i\le \hat{k}} X_i\right| \tag{2.12}$$

based on $X_1, \dots, X_{\hat{k}}$ and

$$T_{n,2} = \frac{1}{s_{n,2}}(n-\hat{k})^{-1/2}\max_{\hat{k}<k\le n}\left|\sum_{\hat{k}<i\le k} X_i - \frac{k-\hat{k}}{n-\hat{k}}\sum_{\hat{k}<i\le n} X_i\right| \tag{2.13}$$

based on $X_{\hat{k}+1}, \dots, X_n$. In (2.12) and (2.13), $s_{n,1}$ and $s_{n,2}$ are equal to the Bartlett estimator computed, respectively, from $X_1, \dots, X_{\hat{k}}$ and $X_{\hat{k}+1}, \dots, X_n$. Specifically, setting

$$\bar{X}_k = \frac{1}{k}\sum_{1\le i\le k} X_i, \qquad \tilde{X}_k = \frac{1}{n-k}\sum_{k<i\le n} X_i$$

and

$$\omega_j(q) = 1 - \frac{j}{q+1}, \tag{2.14}$$

we have

$$s_{n,1}^2 = \frac{1}{\hat{k}}\sum_{1\le i\le \hat{k}} (X_i - \bar{X}_{\hat{k}})^2$$

$$+2\sum_{1\le j\le q(\hat{k})}\omega_j(q(\hat{k}))\frac{1}{\hat{k}}\sum_{1\le i\le \hat{k}-j}(X_i-\bar{X}_{\hat{k}})(X_{i+j}-\bar{X}_{\hat{k}}), \tag{2.15}$$

$$s_{n,2}^2=\frac{1}{n-\hat{k}}\sum_{\hat{k}<i\le n}(X_i-\tilde{X}_{\hat{k}})^2$$
$$+2\sum_{1\le j\le q(n-\hat{k})}\omega_j(q(n-\hat{k}))\frac{1}{n-\hat{k}}\sum_{\hat{k}<i\le n-j}(X_i-\tilde{X}_{\hat{k}})(X_{i+j}-\tilde{X}_{\hat{k}}). \tag{2.16}$$

The test statistic is defined as

$$M_n=\max\{T_{n,1},T_{n,2}\}. \tag{2.17}$$

We first derive the asymptotic distribution of $M_n$ under $H_0$. We need to impose additional assumptions on the change point-model (2.1): both $k^*$, the time of change and $\Delta$, the size of the change, depend on the sample size $n$ such that

$$k^* = [n\theta] \qquad \text{for some } 0<\theta<1, \tag{2.18}$$

$$n\Delta^2 \to \infty, \tag{2.19}$$

$$\Delta^2|\hat{k}-k^*| = O_P(1). \tag{2.20}$$

Condition (2.20) is known to hold if the observations are uncorrelated and was extended by Bai [4], Proposition 3, to moving averages driven by white noise. It also holds if the process $Y_i$ in the change point model (2.1) is strictly stationary, satisfies the approximation condition (2.6), $\Delta\to 0$ and (2.19) holds; see Theorem 4.1.4 in [11]. [There is a misprint in that theorem and $\gamma=0$, which corresponds to our statistic $T_n$, should be included in part (i). The tail condition (4.1.9) in [11] is not needed because it is used only for $\gamma>0$.] Since the squares of ARCH($\infty$) processes satisfy (2.6) (see Theorem 2.1 in [16]), (2.20) holds for such processes.

We will also often impose the following condition on the bandwidth $q(n)$:

$$q(n)\Delta^2=O(1). \tag{2.21}$$

THEOREM 2.1. *Suppose* $H_0$ *and* (2.18)–(2.21) *hold. Suppose* $q(n)$ *is nondecreasing and satisfies*

$$\sup_{k\ge 0}\frac{q(2^{k+1})}{q(2^k)}<\infty, \tag{2.22}$$

$$q(n)\to\infty \quad \textit{and} \quad q(n)(\log n)^4=O(n). \tag{2.23}$$

*Then*

$$(T_{n,1},T_{n,2}) \xrightarrow{d} \Big(\sup_{0\le t\le 1}|B^{(1)}(t)|, \sup_{0\le t\le 1}|B^{(2)}(t)|\Big),$$

*where* $B^{(1)}$ *and* $B^{(2)}$ *are independent Brownian bridges.*

Theorem 2.1 is proved in Appendix B.

COROLLARY 2.1. *Under the assumptions of Theorem 2.1, we have*

$$M_n \xrightarrow{d} \max\Big\{\sup_{0\le t\le 1}|B^{(1)}(t)|, \sup_{0\le t\le 1}|B^{(2)}(t)|\Big\}.$$

Since the distribution function of $\sup_{0\le t\le 1}|B(t)|$ is known (cf. Section 1.5 of [12]), the limit distribution in Corollary 2.1 can be computed explicitly.

In order to describe the asymptotic behavior of the vector $(T_{n,1},T_{n,2})$ if the observations $X_i$ are long-range dependent, we define

$$B_H(t) = W_H(t) - tW_H(1)$$

and

$$\xi = \inf\Big\{t\ge 0 : |B_H(t)| = \sup_{0\le s\le 1}|B_H(s)|\Big\}. \tag{2.24}$$

THEOREM 2.2. *Suppose* $\mathrm{H}_A$ *holds. Assume* $q(n)$ *is nondecreasing, satisfies* (2.22) *and*

$$q(n)\to\infty \quad \text{and} \quad q(n) = O(n(\log n)^{-7/(4-4H)}). \tag{2.25}$$

*Then, the sequence of random vectors*

$$\bigg[\bigg(\frac{q(\hat{k})}{n}\bigg)^{H-1/2} T_{n,1}, \bigg(\frac{q(n-\hat{k})}{n}\bigg)^{H-1/2} T_{n,2}\bigg]$$

*converges in distribution to the random vector*

$$\bigg[\frac{1}{\sqrt{\xi}}\sup_{0\le t\le \xi}\bigg|W_H(t) - \frac{t}{\xi}W_H(\xi)\bigg|,$$
$$\frac{1}{\sqrt{1-\xi}}\sup_{\xi\le t\le 1}\bigg|(W_H(t)-W_H(\xi)) - \frac{t-\xi}{1-\xi}(W_H(1)-W_H(\xi))\bigg|\bigg].$$

Theorem 2.2 is proved in Appendix C.

Theorem 2.2 implies that $T_{n,1}$ and $T_{n,2}$ tend to infinity in probability. Consequently, the test statistic $M_n$ tends to infinity in probability under $\mathrm{H}_A$.

**3. Discussion and examples.** One of the most often used statistics for testing the null hypothesis $\Delta = 0$ in the change-point model (2.1) is the CUSUM statistic

$$T_n = \frac{1}{n^{1/2} s_n} \max_{1 \le k \le n} \left| \sum_{1 \le i \le k} X_i - \frac{k}{n} \sum_{1 \le i \le n} X_i \right|, \tag{3.1}$$

where $s_n^2$ is a suitable estimator of the variance of the sample mean of the $X_i$. If the $Y_i$ in (2.1) are independent identically distributed, $s_n^2$ can be taken to be the sample variance. In this paper we allow the $Y_i$ to be dependent and consider the estimator

$$s_n^2 = \hat{\gamma}_0 + 2 \sum_{1 \le j \le q(n)} \omega_j(q(n)) \hat{\gamma}_j, \tag{3.2}$$

where

$$\hat{\gamma}_j = \frac{1}{n} \sum_{1 \le i \le n-j} (X_i - \bar{X}_n)(X_{i+j} - \bar{X}_n) \tag{3.3}$$

are the sample autocovariances and $\omega_j(q)$ are the Bartlett weights defined by (2.14).

If the observations are weakly dependent (with no change in the mean), the statistic $T_n$ converges to the supremum of a Brownian bridge. However, $T_n \xrightarrow{P} \infty$ either if there is a shift in mean or if the observations are long-range dependent. The latter case is often referred to as a spurious rejection of the null hypothesis of no change in mean. We formalize these observations in Theorems 3.1, 3.2 and 3.3 which, together with Theorems 2.1 and 2.2, form a theoretical foundation for the multistage testing procedure described later in this section. In Theorems 3.1, 3.2 and 3.3, it suffices to assume that

$$q(n) \to \infty \quad \text{and} \quad q(n)/n \to 0 \qquad \text{as } n \to \infty. \tag{3.4}$$

THEOREM 3.1. *Suppose observations $X_1, \dots, X_n$ follow model (2.1) with $\Delta = 0$. If Assumption 2.1 and (3.4) hold, then*

$$T_n \xrightarrow{d} \sup_{0 \le t \le 1} |B(t)|,$$

*where $\{B(t), 0 \le t \le 1\}$ is a Brownian bridge.*

PROOF. Theorem 3.1(i) in [18] implies that if the observations $X_i$ satisfy $X_i = \mu + Y_i$ with the $Y_i$ satisfying Assumption 2.1 and if (3.4) holds, then

$$s_n \xrightarrow{P} \sigma, \tag{3.5}$$

where $\sigma$ is the asymptotic standard deviation appearing in condition (2.3). □

THEOREM 3.2. *Suppose the observations* $X_1, \dots, X_n$ *follow model* (2.1). *If Assumption* 2.1, (3.4), (2.18)–(2.21) *hold, then* $T_n \xrightarrow{P} \infty$. [*Assumption* (2.19) *implies that* $\Delta \neq 0$.]

Theorem 3.2 is proved in Appendix D.

THEOREM 3.3. *Suppose the sequence* $\{X_k\}$ *satisfies Assumption* 2.2. *If* $q(n)/n \to 0$, *then*

$$\left(\frac{q(n)}{n}\right)^{H-1/2} T_n \xrightarrow{d} \sup_{0 \le t \le 1} |W_H(t) - tW_H(1)|. \tag{3.6}$$

[*Convergence* (3.6) *implies* $T_n \xrightarrow{P} \infty$.]

PROOF. By Theorem 3.1 in [18], if (2.9), (2.10) and (3.4) hold, then

$$q(n)^{1-2H} s_n^2 \xrightarrow{P} c_H^2 = \frac{c_0}{H(2H-1)}. \tag{3.7}$$

The constants $c_H$ and $c_0$ in (3.7) are the same as, respectively, in (2.7) and (2.9). Theorem 3.3 now follows immediately from (2.7) and (3.7). □

In order to focus on essential arguments, we considered in Section 2 a simple testing problem. In some applications, however, the presence of more than one change-point may be suspected. Our test can be extended to a multistage testing procedure which is applicable in situations when there is an upper bound on the number of possible change-points. The latter assumption is often used in change-point analysis; see, for example, [44] and references therein. For example, in time series of daily returns on market indices over a period of ten years, or in temperature series over periods of 300 years, one suspects at most two or three change-points; see Section 3 for a data example. For such time series, the maximum number of change points in mean can typically be readily established by a visual inspection of a time series plot.

Before describing the procedure, we must introduce additional notation. Denote by $T(l,m)$ the CUSUM statistic $T_n$ (3.1) computed from the observations $X_{l+1}, \dots, X_m$ and by $\hat{k}(l,m)$ the change-point estimator (2.11) computed from the same observations. Let $B^{(u)}, u = 1, 2, \dots,$ be independent Brownian bridges. Define the critical value $c(u)$ by

$$P\left(\max\left\{\sup_{0 \le t \le 1} |B^{(1)}(t)|, \dots, \sup_{0 \le t \le 1} |B^{(u)}(t)|\right\} > c(u)\right) = \alpha.$$

As mentioned earlier, the distribution of $\sup_{0\le t\le 1}|B^{(1)}(t)|$ is known and is tabulated in [27], so $c(u)$ can be found directly from

$$P\Big(\sup_{0\le t\le 1}|B^{(1)}(t)|\le c(u)\Big)=(1-\alpha)^{1/u}.$$

The procedure we recommend is based on the binary segmentation method of [43]. To focus attention, suppose there can be at most two changes in mean, that is, we want to determine if the observations are weakly dependent with none, one or two changes in the mean or whether they contain a long-range dependent stretch of data. If $T_n=T(0,n)\le c(1)$, the observations are weakly dependent. If $T_n>c(1)$, we compute $\hat{k}_1:=\hat{k}(0,n)$ and $\hat{M}_2=\max[T(0,\hat{k}_1),T(\hat{k}_1,n)]$. If $\hat{M}_2\le c(2)$, the observations are weakly dependent with one change-point. If $\hat{M}_2>c(2)$, we compare $T(0,\hat{k}_1)$ and $T(\hat{k}_1,n)$. Suppose that $T(0,\hat{k}_1)<T(\hat{k}_1,n)$. We then compute $\hat{k}_2=\hat{k}(\hat{k}_1,n)$ and

$$\hat{M}_3=\max[T(0,\hat{k}_1),T(\hat{k}_1,\hat{k}_2),T(\hat{k}_2,n)].$$

Extending Theorem 2.1 to the case of exactly two changes, we have

$$\hat{M}_3\xrightarrow{d}\max\Big\{\sup_{0\le t\le 1}|B^{(1)}(t)|,\sup_{0\le t\le 1}|B^{(2)}(t)|,\sup_{0\le t\le 1}|B^{(3)}(t)|\Big\}.$$

Thus, if $\hat{M}_3\le c(3)$, the observations are weakly dependent with two change-points. If $\hat{M}_3>c(3)$, the observations contain a long-range dependent stretch of data.

Before concluding this section with a data example, we list several time series models which satisfy Assumptions 2.1 or 2.2. References to the proofs can be found in [16].

EXAMPLE 3.1. The linear process

$$X_k=\sum_j a_j\varepsilon_{k-j}, \tag{3.8}$$

where $\varepsilon_j$ are independent identically distributed random variables with finite fourth moment and zero mean, satisfies Assumption 2.1 if $\sum_j|a_j|<\infty$. In particular, ARMA processes whose autoregressive polynomial has no zeros on the unit circle satisfy Assumption 2.1.

If, on the other hand, $a_j\sim cj^{d-1}$ for some $0<d<1/2$, then the $X_k$ (3.8) satisfy Assumption 2.2 with $H=d+1/2$. In particular, fractional ARIMA processes whose autoregressive polynomial has no zeros on the unit circle satisfy Assumption 2.2.

EXAMPLE 3.2. Consider the process $\{\eta_k\}$ satisfying

$$\eta_k = \rho_k \xi_k, \qquad \rho_k = \sum_{j \geq 1} c_j \eta_{k-j}, \tag{3.9}$$

where $a > 0, c_j \geq 0$ and the $\xi_k$ are independent identically distributed nonnegative random variables with finite fourth moment. The $\eta_k$ should be viewed as the squares of an ARCH process. If

$$[E\xi_0^4]^{1/4} \sum_{j \geq 1} c_j < 1, \tag{3.10}$$

then the sequence $Y_k = \eta_k - E\eta_k$ satisfies Assumption 2.1.

As a more specific example, consider $Y_k = r_k^2 - Er_k^2$, where the $r_k$ follow a GARCH$(p, q)$ model,

$$r_k = \sigma_k \varepsilon_k, \qquad \sigma_k^2 = \omega + \sum_{1 \leq i \leq p} \alpha_i r_{k-i}^2 + \sum_{1 \leq j \leq q} \beta_j \sigma_{k-j}^2. \tag{3.11}$$

Then, under regularity conditions derived in [7], the $c_i$ are defined by

$$\sum_{j \geq 1} c_j z^j = \frac{\sum_{1 \leq i \leq p} \alpha_i z^i}{1 - \sum_{1 \leq j \leq q} \beta_j z^j}, \qquad |z| \leq 1,$$

and $\rho_k = \sigma_k^2, \xi_k = \varepsilon_k^2$.

EXAMPLE 3.3. The $r_k$ are said to follow a LARCH (Linear ARCH) model if

$$r_k = \sigma_k \varepsilon_k, \qquad \sigma_k = a + \sum_{j \geq 1} b_j r_{k-j}, \tag{3.12}$$

where $a \neq 0$, the $b_j$ are real coefficients (not necessarily nonnegative) and the $\varepsilon_k$ are independent identically distributed with zero mean and finite fourth moment. If $b_j \sim c j^{d-1}$ for some $0 < d < 1/2$ and

$$L[E\varepsilon_0^4]^{1/2} \sum_{j \geq 1} b_j^2 < 1,$$

where $L = 7$ if the $\varepsilon_k$ are Gaussian and $L = 11$ in general, then $Y_k = r_k^2$ satisfy conditions (2.7) and (2.9) of Assumption 2.2. Conditions for (2.10) to hold have not been established yet. The LARCH model was studied by Robinson [39], Giraitis et al. [17, 19, 20] and Berkes and Horváth [6], among others.

We conclude this section with an illustration of how our procedure can be applied in practice. Figure 1 shows daily returns of the Dow Jones Industrial Average from January 1, 1992 to December 31, 1999 and a simulated

LARCH process with $H = 0.85$ of the same length ($n = 2021$). The corresponding columns show the sample autocorrelation functions and smoothed periodograms of the squares of the two series in the top row. The volatility (variance) of the Dow Jones series appears to have a change point somewhere in the middle of series, but given that we observe only a finite realization, this change-point might be spurious and the observed change in variance might be explained as a *persistent* increase in volatility characteristic of a long memory process. That this might well be the case is reinforced by the examination of the plot of the simulated LARCH series which exhibits markedly higher variability in the first 1/3 of the realization, even though the plot shows a realization of a strictly and fourth-order stationary process. The left column in Figure 1 shows that the autocorrelation function of the squared Dow Jones returns does not decay to zero in a fashion typical of a short memory process and the smoothed periodogram (on a log–log scale) exhibits a clear positive slope. In fact, a periodogram-based semipara-metric estimate of $H$ based on the automatic bandwith selection procedure proposed by Lobato and Robinson [32] yields the estimate $\hat{H} = 0.842991$.

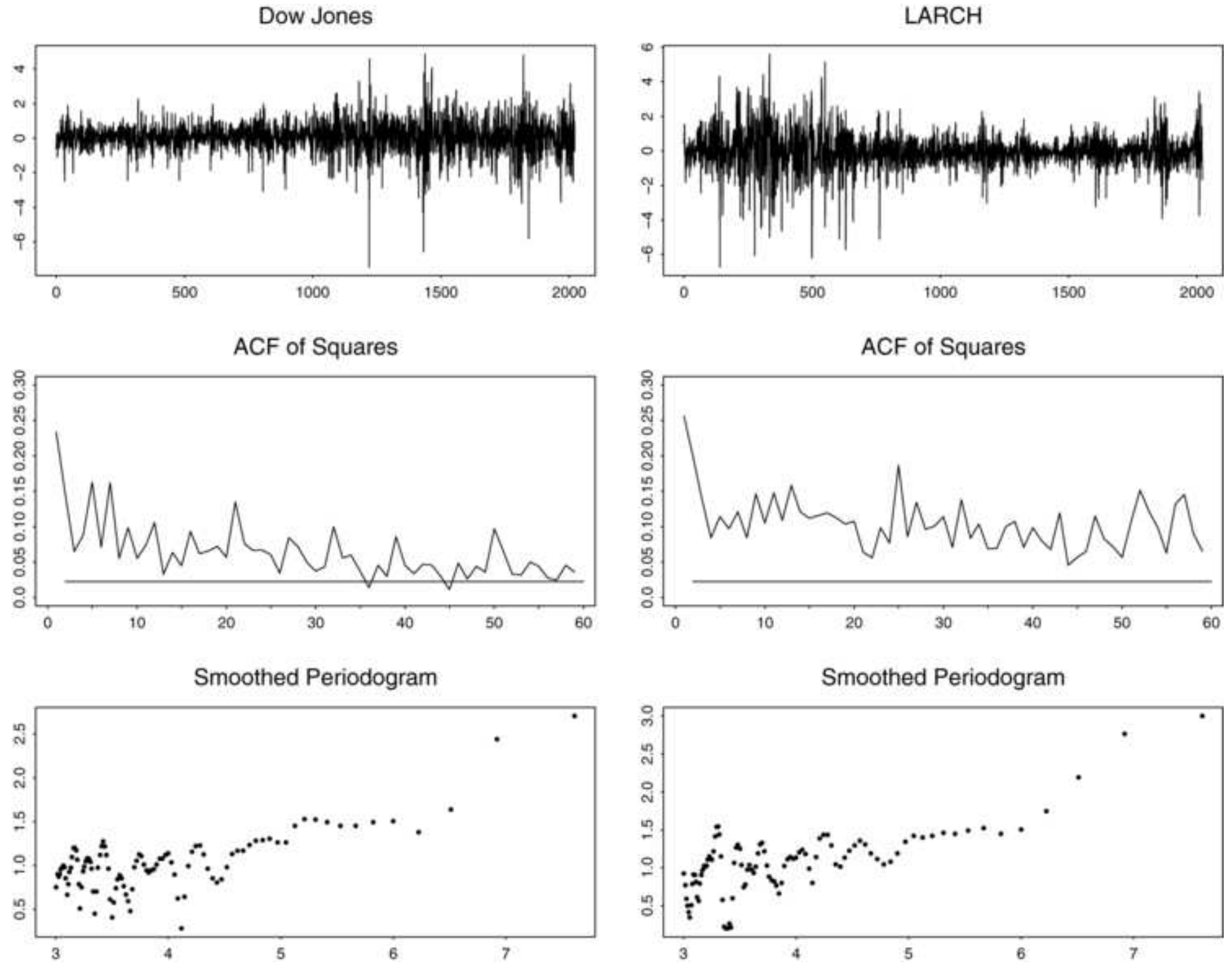

FIG. 1. *Daily returns of the Dow Jones Industrial Average and a simulated LARCH process with $H = 0.85$ together with the autocorrelation functions and smoothed periodograms at low frequencies of the squared observations.*

| | $\boldsymbol{\omega}$ | $\boldsymbol{\alpha_1}$ | $\boldsymbol{\beta_1}$ |
|---|---|---|---|
| Before $\hat{k}$ | 0.02461474 | 0.06404848 | 0.87864088 |
| After $\hat{k}$ | 0.09540076 | 0.09734341 | 0.83945713 |

For the Dow Jones returns, we therefore wish to test the null hypothesis of exactly one change in the variance of the observations against the alternative that the squared observations are a realization of a long-range dependent process. Assuming that the mean of the returns is zero (we subtracted the sample mean of 0.05829482 before conducting further analysis), this testing problem is thus identical with the basic testing problem formulated in Section 2, with the $X_i$ being equal to the squared returns.

In order to perform the test, we need to choose the bandwidth function $q(\cdot)$. We performed our calculations in `Splus` and used the function `acf` to obtain sample autocovariances. By default, this function returns the first $10 \log_{10}(n)$ sample autocovariances for a time series of length $n$. We found, however, that, for the nonlinear return data, more autocovariances must be used to capture the dependence structure, so we increased the maximum lag up to which the autocovariances are computed by 50%. Thus, in the following, we report the results based on

$$q(n) = 15 \log_{10}(n).$$

The value of $M_n$ is 1.341153, which lies below the 10% asymptotic critical value of 1.36 (the 5% and 1% critical values are, resp., 1.48 and 1.72). We are thus unable to reject the null hypothesis of a change-point in the level of the squared returns.

To validate the above conclusion, we need to assess the empirical size and power of the test. To assess the size, we divided the data into two parts: before and after the estimated change point $\hat{k} = 1061$ and fitted the GARCH(1, 1) model to each stretch of data. We obtained the following parameters [see (3.11)]:
Using $Er_k^2 = \omega/(1 - \alpha_1 - \beta_1)$, the implied change in the variance (level of the $r_k^2$) is 1.080022. In fact, variances implied by the GARCH(1, 1) models before and after $\hat{k}$ are very close to the corresponding sample variances whose difference is 1.040886.

We simulated one thousand replications of the above change point model and on each of them we computed the value of $M_n$. Table 1 reports the percentage of rejections of the null hypothesis. At the nominal confidence level of 10%, the percentage of rejections is slightly over 10%, suggesting that accepting the null hypothesis based on the value of $M_n = 1.34115$ was not due to type II error.

To assess the power, we simulated one thousand replications of the LARCH process (3.12) with $d = 0.35$ and the $b_j$ computed according to the recursion $b_j = [b_{j-1}(j+d)]/(j+1)$, with $b_0 = 0.25$ and $a = 0.03$. These parameter values ensure that the process is fourth-order stationary and were chosen by experimentation to make the realizations similar to the Dow Jones returns, with a typical realization shown in Figure 1. Table 1 shows that the test is able to detect the alternative at the nominal 10% level with probability of over 30%. For this particular alternative, the power is not very high. This can be explained by the fact that the realizations of a LARCH process with the parameters chosen above and for the sample size of $n = 2021$ often exhibit two periods of different variability which can by separated by the change-point estimator $\hat{k}$. The intensity of long-range dependence in each of the two subsamples is "underestimated," yielding small values of $M_n$. However, even though the alternative is "very close" to the null, the test has nontrivial power.

The above illustration is not meant as a guide for practitioners, but merely points out the potential of the test.

## APPENDIX A

**Almost sure convergence of the Bartlett estimator.** For ease of reference, we present here the result on the almost sure asymptotics for the estimator $s_n^2$, which we appeal to in the folllowing. Its proof is given in [8].

THEOREM A.1. *Suppose $\{Y_k\}$ is a fourth-order stationary sequence with $EY_i = 0$ and $\gamma_j = \mathrm{Cov}(Y_0, Y_j)$. Consider the variance estimator*

$$s_n^2 = \hat{\gamma}_0 + 2 \sum_{1 \le j \le q(n)} \omega_j(q(n)) \hat{\gamma}_j, \tag{A.1}$$

*where $\hat{\gamma}_j$ are the sample autocovariances and $\omega_j(q)$ are the Bartlett weights defined respectively by (3.3) and (2.14).*

*Suppose the sequence $q(n)$ is nondecreasing and*

$$\sup_{k \ge 0} \frac{q(2^{k+1})}{q(2^k)} < \infty. \tag{A.2}$$

TABLE 1
*Empirical size and power of the asymptotic test based on the statistic $M_n$*

| **Nominal level (in %)** | **10.0** | **5.0** | **1.0** |
|---|---|---|---|
| Empirical size | 13.4 | 6.5 | 0.8 |
| Empirical power | 32.5 | 20.0 | 5.0 |

(i) *Suppose, in addition, that conditions* (2.4) *and* (2.5) *hold and*

$$q(n) \to \infty \quad \textit{and} \quad q(n)(\log n)^4 = O(n). \tag{A.3}$$

*Then*

$$s_n^2 \to \sigma^2 := \sum_{j=-\infty}^{\infty} \gamma_j \qquad \textit{a.s.} \tag{A.4}$$

(ii) *Assume*

$$\tfrac{1}{2} < H < 1 \tag{A.5}$$

*and*

$$\gamma_k \sim c_0 k^{2H-2} \tag{A.6}$$

*for some* $c_0 > 0$. *Assume also that*

$$q(n) \to \infty \quad \textit{and} \quad q(n) = O(n(\log n)^{-7/(4-4H)}) \tag{A.7}$$

*and*

$$\sup_{|h| \le q(n)} \sum_{-n \le r,s \le n} |\kappa(h,r,s)| = O(n^{2H-1}). \tag{A.8}$$

*Then*

$$q(n)^{1-2H} s_n^2 \to c_H^2 = \frac{c_0}{H(2H-1)} \qquad \textit{a.s.} \tag{A.9}$$

REMARK A.1. By the fourth-order stationarity of the $X_i$, all bounds in the proof of Theorem A.1 remain valid if the random variables $X_1, \dots, X_n$ are replaced by $X_{k+1}, \dots, X_n$, $n$ is replaced by $n-k$ and $q(n)$ is replaced by $q(n-k)$. Therefore, on denoting $\tilde{X}_k = \frac{1}{n-k} \sum_{k<i\le n} X_i$ and

$$\begin{aligned} s_{k,n}^2 &= \frac{1}{n-k} \sum_{k<i\le n} (X_i - \tilde{X}_k)^2 \\ &\quad + 2 \sum_{1\le j\le q(n-k)} \omega_j(q(n-k)) \frac{1}{n-k} \sum_{k<i\le n} (X_i - \tilde{X}_k)(X_{i+j} - \tilde{X}_k), \end{aligned} \tag{A.10}$$

under the assumptions of part (i) of Theorem A.1,

$$s_{k,n}^2 \overset{\text{a.s.}}{\to} \sigma^2 \qquad \text{as } n-k \to \infty, \tag{A.11}$$

and under the assumptions of part (ii) of Theorem A.1,

$$[q(n-k)]^{1-2H} s_{k,n}^2 \overset{\text{a.s.}}{\to} c_H^2 \qquad \text{as } n-k \to \infty. \tag{A.12}$$

Relations (A.11) and (A.12) are used, respectively, in the proofs of Lemmas B.3 and C.2.

## APPENDIX B

**Proof of Theorem 2.1.** Theorem 2.1 will follow immediately from Lemmas B.1, B.2 and B.3 which are stated and proved below.

In this section we assume that the observations follow the change-point model (2.1) and that (2.18) holds.

We will extensively use the relation

$$\left|\frac{\hat{k}}{n} - \theta\right| = o_P(1), \tag{B.1}$$

which follows from assumptions (2.19) and (2.20).

Lemma B.1. *If (2.19) and (2.20) hold, then*

$$\begin{aligned} &n^{-1/2} \max_{1\le k\le \hat{k}} \left| \sum_{1\le i\le k} X_i - \frac{k}{\hat{k}} \sum_{1\le i\le \hat{k}} X_i \right| \\ &\quad = n^{-1/2} \max_{1\le k\le \hat{k}} \left| \sum_{1\le i\le k} Y_i - \frac{k}{\hat{k}} \sum_{1\le i\le \hat{k}} Y_i \right| + o_P(1) \end{aligned} \tag{B.2}$$

*and*

$$\begin{aligned} &n^{-1/2} \max_{\hat{k}< k\le n} \left| \sum_{\hat{k}< i\le k} X_i - \frac{k-\hat{k}}{n-\hat{k}} \sum_{\hat{k}< i\le n} X_i \right| \\ &\quad = n^{-1/2} \max_{\hat{k}< k\le n} \left| \sum_{\hat{k}< i\le k} Y_i - \frac{k-\hat{k}}{n-\hat{k}} \sum_{\hat{k}< i\le n} Y_i \right| + o_P(1). \end{aligned} \tag{B.3}$$

Proof. We can assume that $\mu = 0$. Since the verification of (B.3) is very similar to that of (B.2), we present only the proof of (B.2).

If $\hat{k} \le k^*$,

$$\sum_{1\le i\le k} X_i - \frac{k}{\hat{k}} \sum_{1\le i\le \hat{k}} X_i = \sum_{1\le i\le k} Y_i - \frac{k}{\hat{k}} \sum_{1\le i\le \hat{k}} Y_i$$

for all $1 \le k \le \hat{k}$, so (B.2) holds trivially. If $k^* < \hat{k}$, then

$$\frac{k}{\hat{k}} \sum_{1\le i\le \hat{k}} X_i = \frac{k}{\hat{k}} \sum_{1\le i\le \hat{k}} Y_i + \Delta \frac{k}{\hat{k}} (\hat{k} - k^*)$$

and

$$\sum_{1\le i\le k} X_i = \begin{cases} \displaystyle\sum_{1\le i\le k} Y_i, & \text{if } 1 \le k \le k^*, \\ \displaystyle\sum_{1\le i\le k} Y_i + (k-k^*)\Delta, & \text{if } k^* < k \le \hat{k}. \end{cases}$$

Hence,

$$\left| n^{-1/2} \max_{1\le k\le \hat{k}} \left| \sum_{1\le i\le k} X_i - \frac{k}{\hat{k}} \sum_{1\le i\le \hat{k}} X_i \right| - n^{-1/2} \max_{1\le k\le \hat{k}} \left| \sum_{1\le i\le k} Y_i - \frac{k}{\hat{k}} \sum_{1\le i\le \hat{k}} Y_i \right| \right|$$

$$\le 2n^{-1/2}\Delta|\hat{k} - k^*| = \frac{2\Delta^2|\hat{k} - k^*|}{\Delta n^{1/2}},$$

so (B.2) follows from assumptions (2.20) and (2.19). □

LEMMA B.2. *If (2.3) and (B.1) hold, then the sequence of random vectors*

$$\left( \hat{k}^{-1/2} \max_{1\le k\le \hat{k}} \left| \sum_{1\le i\le k} Y_i - \frac{k}{\hat{k}} \sum_{1\le i\le \hat{k}} Y_i \right|, \right.$$

$$\left. (n-\hat{k})^{-1/2} \max_{\hat{k}<k\le n} \left| \sum_{\hat{k}<i\le k} Y_i - \frac{k-\hat{k}}{n-\hat{k}} \sum_{\hat{k}<i\le n} Y_i \right| \right)$$

*converges in distribution to the random vector*

$$\left( \sup_{0\le t\le 1} |B^{(1)}(t)|, \sup_{0\le t\le 1} |B^{(2)}(t)| \right),$$

*where $B^{(1)}$ and $B^{(2)}$ are independent Brownian bridges.*

PROOF. By (2.6),

$$\max_{1\le k\le \hat{k}} \left| \left[ \sum_{1\le i\le k} Y_i - \frac{k}{\hat{k}} \sum_{1\le i\le \hat{k}} Y_i \right] - \sigma n^{1/2} \left[ W_n\left(\frac{k}{n}\right) - \frac{k}{\hat{k}} W_n\left(\frac{\hat{k}}{n}\right) \right] \right|$$

$$= o_P(n^{1/2}). \tag{B.4}$$

Using (B.1) and the continuity of the Wiener process, we get

$$|W_n(\hat{k}/n) - W_n(\theta)| = o_P(1). \tag{B.5}$$

Hence,

$$\max_{1\le k\le \hat{k}} \left| W_n\left(\frac{k}{n}\right) - \frac{k}{\hat{k}} W_n\left(\frac{\hat{k}}{n}\right) \right| = \sup_{0\le t\le \theta} \left| W_n(t) - \frac{t}{\theta} W_n(\theta) \right| + o_P(1), \tag{B.6}$$

and consequently, by (B.4),

$$\hat{k}^{-1/2} \max_{1\le k\le \hat{k}} \left| \sum_{1\le i\le k} Y_i - \frac{k}{\hat{k}} \sum_{1\le i\le \hat{k}} Y_i \right|$$

$$= \frac{\sigma}{\theta^{1/2}} \sup_{0\le t\le \theta} \left| W_n(t) - \frac{t}{\theta} W_n(\theta) \right| + o_P(1). \tag{B.7}$$

Similar arguments give

$$
\begin{aligned}
&(n-\hat{k})^{-1/2} \max_{\hat{k}<k\le n}\left|\sum_{\hat{k}<i\le k} Y_i - \frac{k-\hat{k}}{n-\hat{k}}\sum_{\hat{k}<i\le n} Y_i\right| \\
&\quad = \frac{\sigma}{(1-\theta)^{1/2}} \sup_{\theta\le t\le 1}\left|(W_n(t)-W_n(\theta)) \right. \\
&\qquad\qquad \left. - \frac{t-\theta}{1-\theta}(W_n(1)-W_n(\theta))\right| + o_P(1).
\end{aligned}
\tag{B.8}
$$

Since $\theta^{-1/2}W_n(\theta t), 0\le t\le 1$, is a Wiener process,

$$
\frac{1}{\theta^{1/2}} \sup_{0\le t\le\theta}\left|W_n(t)-\frac{t}{\theta}W_n(\theta)\right| \stackrel{d}{=} \sup_{0\le t\le 1}|B^{(1)}(t)|,
\tag{B.9}
$$

where $B^{(1)}$ is a Brownian bridge. Similarly, there is a Wiener process $W(t), 0\le t\le 1$, such that

$$
\begin{aligned}
&\frac{1}{(1-\theta)^{1/2}} \sup_{\theta\le t\le 1}\left|(W_n(t)-W_n(\theta)) - \frac{t-\theta}{1-\theta}(W_n(1)-W_n(\theta))\right| \\
&\quad \stackrel{d}{=} \frac{t-\theta}{1-\theta} \sup_{\theta\le t\le 1}\left|W(t-\theta)-\frac{t-\theta}{1-\theta}W(1-\theta)\right| \\
&\quad = \frac{t-\theta}{1-\theta} \sup_{0\le t\le 1-\theta}\left|W(t)-\frac{t}{1-\theta}W(1-\theta)\right| \stackrel{d}{=} \sup_{0\le t\le 1}|B^{(2)}(t)|,
\end{aligned}
\tag{B.10}
$$

where $B^{(2)}$ is another Brownian bridge. The claim thus follows by combining (B.7), (B.8) and (B.9), (B.10) and using the independence of the increments of a Wiener process. □

LEMMA B.3. *Suppose Assumption* 2.1, (2.19), (2.21), (2.22), (2.23) *and* (B.1) *hold. Then*

$$
s_{n,1} \xrightarrow{P} \sigma \quad \textit{and} \quad s_{n,2} \xrightarrow{P} \sigma.
$$

PROOF. Following the proof of Proposition D.1, we get

$$
s_{n,1}^2 = \sum_{1\le m\le 5}\left[\hat{\gamma}_{0m,1} + 2\sum_{1\le j\le q(\hat{k})} \omega_j(q(\hat{k}))\hat{\gamma}_{jm,1}\right] =: \sum_{1\le m\le 5} s_{nm,1}^2,
$$

where

$$
\hat{\gamma}_{j1,1} = \frac{1}{\hat{k}} \sum_{1\le i\le \hat{k}-j} (Y_i - \bar{Y}_{\hat{k}})(Y_{i+j}-\bar{Y}_{\hat{k}}),
$$

$$\hat{\gamma}_{j2,1} = \frac{1}{\hat{k}}\left[(k^* - j)\left(\frac{\hat{k}-k^*}{n}\Delta\right)^2 - j\frac{\hat{k}-k^*}{\hat{k}}\frac{k^*}{n}\Delta^2 + (\hat{k}-j-k^*)\left(\frac{k^*}{\hat{k}}\Delta\right)^2\right],$$

$$\hat{\gamma}_{j3,1} = -\frac{1}{\hat{k}}\sum_{1\le i\le k^*-j}[(Y_i - \bar{Y}_{\hat{k}}) + (Y_{i+j} - \bar{Y}_{\hat{k}})]\frac{n-k^*}{n}\Delta,$$

$$\hat{\gamma}_{j4,1} = \frac{1}{\hat{k}}\sum_{k^*-j<i\le k^*}\left[(Y_i - \bar{Y}_{\hat{k}})\frac{k^*}{\hat{k}} - (Y_{i+j} - \bar{Y}_{\hat{k}})\frac{\hat{k}-k^*}{\hat{k}}\right]\Delta$$

and

$$\hat{\gamma}_{j5,1} = \frac{1}{\hat{k}}\sum_{k^*<i\le\hat{k}-j}[(Y_i - \bar{Y}_{\hat{k}}) + (Y_{i+j} - \bar{Y}_{\hat{k}})]\frac{k^*}{\hat{k}}\Delta.$$

Since $s_n^2 \to \sigma^2$ a.s. by part (i) of Theorem A.1 and $\hat{k} \xrightarrow{P} \infty$ by (B.1), Theorem 7.1.1(c) on page 252 of [12] yields that

$$s_{n1,1}^2 \xrightarrow{P} \sigma^2. \tag{B.11}$$

Next we show that

$$s_{nm,1}^2 = o_P(1) \qquad \text{for } m = 2,3,4,5. \tag{B.12}$$

As we have seen in the proof of Proposition D.1,

$$s_{nm,1}^2 = O_P(\hat{k}^{-1/2}q(\hat{k})\Delta) = O_P\left(\frac{q(\hat{k})\Delta^2}{\hat{k}^{1/2}\Delta}\right) = O_P\left(\left(\frac{n}{\hat{k}}\right)^{1/2}\frac{q(\hat{k})\Delta^2}{n^{1/2}\Delta}\right) = o_P(1),$$

proving (B.12).

To prove $s_{n2}^2 \xrightarrow{P} \sigma^2$, we can apply the same argument, upon observing that by Remark A.1, for all $0 < r < 1$, we have

$$\max_{rn\le k\le n}|s_{k,n}^2 - \sigma^2| \to 0 \qquad \text{a.s.},$$

where $s_{k,n}^2$ is defined in (A.10). □

## APPENDIX C

**Proof of Theorem 2.2.** Theorem 2.2 will follow directly from Lemmas C.1 and C.2 below. We can assume that $\mu = 0$.

Let

$$Z_{n1}(t) = \frac{1}{n^H}\max_{1\le k\le nt}\left|\sum_{1\le i\le k}X_i - \frac{k}{nt}\sum_{1\le i\le nt}X_i\right|$$

and

$$Z_{n2}(t) = \frac{1}{n^H}\max_{nt<k\le n}\left|\sum_{nt\le i\le k}X_i - \frac{k-nt}{n-nt}\sum_{nt<i\le n}X_i\right|.$$

Similarly, let

$$Z_1(t) = c_H \sup_{0 \le s \le t} \left| W_H(s) - \frac{s}{t} W_H(t) \right|$$

and

$$Z_2(t) = c_H \sup_{t < s \le 1} \left| (W_H(s) - W_H(t)) - \frac{s-t}{1-t}(W_H(1) - W_H(t)) \right|,$$

where $W_H$ is defined in Assumption 2.2.

LEMMA C.1. *Suppose that* (2.7) *and* (2.8) *hold. Then*

$$(\hat{k}/n, Z_{n1}(t), Z_{n2}(t)) \xrightarrow{d} (\xi, Z_1(t), Z_2(t)), \tag{C.1}$$

*where* $\xi$ *is defined by* (2.24). *The vectors in* (C.1) *take values in* $(0,1) \times D[0,1] \times D[0,1]$.

PROOF. The vector $(\hat{k}/n, Z_{n1}(t), Z_{n2}(t))$ is a continuous mapping of $\{n^{-H} \sum_{1 \le k \le nt} X_k, 0 \le t \le 1\}$. The same mapping transforms $W_H(t)$ into $(\xi, Z_1(t), Z_2(t))$. Hence, the statement of the lemma follows from the continuous mapping theorem. $\square$

LEMMA C.2. *We assume that the conditions of Theorem* 2.2 *are satisfied. Then*

$$[q(\hat{k})]^{1-2H} s_{n,1}^2 \xrightarrow{P} c_H^2 \tag{C.2}$$

*and*

$$[q(n-\hat{k})]^{1-2H} s_{n,2}^2 \xrightarrow{P} c_H^2. \tag{C.3}$$

PROOF. We first verify (C.2). By part (ii) of Theorem A.1, for any $0 < r < 1$,

$$\sup_{k \ge rn} |[q(k)]^{1-2H} s_k^2 - c_H^2| \xrightarrow{\text{a.s.}} 0.$$

For any $0 < r < 1$ which is a continuity point of the distribution function of $\xi$ we have

$$\begin{aligned} &\limsup_{n\to\infty} P[|[q(\hat{k})]^{1-2H} s_{\hat{k}}^2 - c_H^2| > \varepsilon] \\ &\qquad \le \limsup_{n\to\infty} P[\hat{k}/n \le r] + \limsup_{n\to\infty} P\left[\sup_{k \ge rn} |[q(k)]^{1-2H} s_k^2 - c_H^2| > \varepsilon\right] \\ &\qquad = P(\xi \le r). \end{aligned}$$

Since $P(\xi \le r) \to 0$ as $r \to 0$, (C.2) follows.

To prove (C.3), note that by (A.12),

$$\sup_{k \le (1-r)n} |[q(n-k)]^{1-2H} s_{k,n}^2 - c_H^2| \overset{\text{a.s.}}{\to} 0.$$

Relation (C.3) is then established using a lim sup argument as above and the fact that $P(\xi > r) \to 0$ as $r \to 1$. □

## APPENDIX D

**Proof of Theorem 3.2.** Observe that by Assumption 2.1 and (2.18),

$$\begin{aligned}
&\frac{1}{n^{1/2}} \max_{1 \le k \le n} \left| \sum_{1 \le i \le k} X_i - \frac{k}{n} \sum_{1 \le i \le n} X_i \right| \\
&\quad \ge \frac{1}{n^{1/2}} \left| \sum_{1 \le i \le k^*} X_i - \frac{k^*}{n} \sum_{1 \le i \le n} X_i \right| \\
&\quad = \frac{1}{n^{1/2}} \left| \sum_{1 \le i \le k^*} Y_i - \frac{k^*}{n} \sum_{1 \le i \le n} Y_i - \frac{k^*(n-k^*)}{n} \Delta \right| \\
&\quad \ge \frac{1}{2} n^{1/2} \theta(1-\theta)|\Delta| - O_P(1),
\end{aligned}$$

as $n \to \infty$. Hence, it suffices to show that

$$\frac{n^{1/2}|\Delta|}{s_n} \overset{P}{\to} \infty,$$

which, in view of (2.19), will follow if we show that $s_n = O_P(1)$, which is verified in the following proposition.

Proposition D.1. *Suppose model (2.1) is valid. Consider the estimator $s_n^2$ defined by (3.2). Suppose Assumption 2.1 and (3.4), (2.18)–(2.21) hold. Then*

$$s_n^2 = O_P(1). \tag{D.1}$$

Proof. Denoting $V_{i,j} = (X_i - \bar{X}_n)(X_{i+j} - \bar{X}_n)$, observe that

$$\begin{aligned}
V_{i,j} &= (Y_i - \bar{Y}_n)(Y_{i+j} - \bar{Y}_n) - (Y_i - \bar{Y}_n)\frac{n-k^*}{n}\Delta \\
&\quad - (Y_{i+j} - \bar{Y}_n)\frac{n-k^*}{n}\Delta + \left(\frac{n-k^*}{n}\Delta\right)^2 \qquad \text{if } 1 \le i \le i+j \le k^*, \\
V_{i,j} &= (Y_i - \bar{Y}_n)(Y_{i+j} - \bar{Y}_n) + (Y_i - \bar{Y}_n)\frac{k^*}{n}\Delta
\end{aligned}$$

$$+ (Y_{i+j} - \bar{Y}_n)\frac{n-k^*}{n}\Delta - \frac{n-k^*}{n}\Delta\frac{k^*}{n}\Delta \qquad \text{if } 1 \le i \le k^* < i+j,$$

$$V_{i,j} = (Y_i - \bar{Y}_n)(Y_{i+j} - \bar{Y}_n) + (Y_i - \bar{Y}_n)\frac{k^*}{n}\Delta$$

$$+ (Y_{i+j} - \bar{Y}_n)\frac{k^*}{n}\Delta + \left(\frac{k^*}{n}\Delta\right)^2 \qquad \text{if } k^* < i \le i+j.$$

Therefore, for any $0 \le j \le q$,

$$\hat{\gamma}_j = \sum_{1 \le m \le 5} \hat{\gamma}_{jm},$$

where

$$\hat{\gamma}_{j1} = \frac{1}{n}\sum_{1 \le i \le n-j}(Y_i - \bar{Y}_n)(Y_{i+j} - \bar{Y}_n),$$

$$\hat{\gamma}_{j2} = \frac{1}{n}\left[(k^* - j)\left(\frac{n-k^*}{n}\Delta\right)^2 - j\frac{n-k^*}{n}\frac{k^*}{n}\Delta^2 + (n-j-k^*)\left(\frac{k^*}{n}\Delta\right)^2\right],$$

$$\hat{\gamma}_{j3} = -\frac{1}{n}\sum_{1 \le i \le k^*-j}[(Y_i - \bar{Y}_n) + (Y_{i+j} - \bar{Y}_n)]\frac{n-k^*}{n}\Delta,$$

$$\hat{\gamma}_{j4} = \frac{1}{n}\sum_{k^*-j < i \le k^*}\left[(Y_i - \bar{Y}_n)\frac{k^*}{n} - (Y_{i+j} - \bar{Y}_n)\frac{n-k^*}{n}\right]\Delta$$

and

$$\hat{\gamma}_{j5} = \frac{1}{n}\sum_{k^* < i \le n-j}[(Y_i - \bar{Y}_n) + (Y_{i+j} - \bar{Y}_n)]\frac{k^*}{n}\Delta.$$

Consequently,

$$s_n^2 = \sum_{1 \le m \le 5}\left[\hat{\gamma}_{0m} + 2\sum_{1 \le j \le q}\omega_j(q)\hat{\gamma}_{jm}\right] =: \sum_{1 \le m \le 5} s_{nm}^2.$$

By (3.5),

$$s_{n1}^2 \xrightarrow{P} \sigma^2. \tag{D.2}$$

Hence, (D.1) will follow if we verify that

$$s_{nm}^2 = O_P(1) \qquad \text{for } m = 2, 3, 4, 5. \tag{D.3}$$

In order to verify (D.3), we will often appeal to the two elementary relations

$$2\sum_{1 \le j \le q}\omega_j(q) \sim q \tag{D.4}$$

and

$$2\sum_{1\le j\le q} j\omega_j(q) \sim \tfrac{1}{3}q^2. \tag{D.5}$$

Relation (D.3) is easy to verify for $m=3,4,5$. By (2.6) and (D.4),

$$s_{nm}^2 = O_P(n^{-1/2}q\Delta) = O_P\left(\frac{q\Delta^2}{n^{1/2}\Delta}\right) = o_P(1)$$

on account of (2.19) and (2.21). It thus remains to establish (D.4) for $m=2$. Since there are three terms in the definition of $\hat{\gamma}_{j2}$, we may write

$$\begin{aligned}\frac{s_{n21}^2}{q\Delta^2} &= \frac{1}{q}\left[\frac{k^*}{n}\left(\frac{n-k^*}{n}\right)^2 + 2\sum_{1\le j\le q}\omega_j(q)\left(\frac{k^*-j}{n}\right)\left(\frac{n-k^*}{n}\right)^2\right]\\ &\sim \frac{1}{q}\left[\theta(1-\theta)^2 + 2\sum_{1\le j\le q}\omega_j(q)\theta(1-\theta)^2\right.\\ &\qquad\qquad \left. - \frac{2}{n}\sum_{1\le j\le q} j\omega_j(q)(1-\theta)^2\right]\\ &\sim \frac{1}{q}\left[(1+q)\theta(1-\theta)^2 - \frac{1}{3}\frac{q^2}{n}(1-\theta)^2\right] \to \theta(1-\theta)^2.\end{aligned} \tag{D.6}$$

Similarly,

$$\frac{s_{n22}^2}{q\Delta^2} = -\frac{1}{q}\frac{n-k^*}{n}\frac{k^*}{n}\frac{2}{n}\sum_{1\le j\le q} j\omega_j(q) \sim -(1-\theta)\theta\frac{q}{3n} \to 0 \tag{D.7}$$

and

$$\begin{aligned}\frac{s_{n23}^2}{q\Delta^2} &= \frac{1}{q}\left[\frac{n-k^*}{n}\left(\frac{k^*}{n}\right)^2 + 2\sum_{1\le j\le q}\omega_j(q)\frac{n-k^*-j}{n}\left(\frac{k^*}{n}\right)^2\right]\\ &\sim \frac{1}{q}\left[(1-\theta)\theta^2 + 2\sum_{1\le j\le q}\omega_j(q)(1-\theta)\theta^2 - \frac{2}{n}\sum_{1\le j\le q} j\omega_j(q)\theta^2\right]\\ &\sim \frac{1}{q}\left[(1+q)(1-\theta)\theta^2 - \frac{\theta^2}{3}\frac{q^2}{n}\right] \to (1-\theta)\theta^2.\end{aligned} \tag{D.8}$$

Putting together relations (D.6), (D.7) and (D.8), we obtain (D.3) for $m=2$. This completes the proof of Proposition D.1. □

REMARK D.1. Proposition D.1 and, therefore Theorem 3.2, remain valid if the Bartlett weights (2.14) are replaced by any weights satisfying (D.4) and $\sum_{1\le j\le q} j\omega_j(q) = O(q^2)$ in addition to the following conditions which are needed for (D.2) to hold: $\omega_j(q) = 0$ for $|j| > q$, $0 \le \omega_j(q) \le 1$, and

$$\lim_{q\to\infty} \omega_j(q) = 1 \qquad \text{for each } j; \tag{D.9}$$

see Remark 1.2 in [8].

**Acknowledgments.** The paper has benefited from the comments of the three referees and constructive and detailed advice of the Associate Editor. We are especially indebted to the first referee for an exceptionally careful reading of the original version of this paper and an accurate and instructive report. The numerical work in Section 3 was done by Aonan Zhang.

I. Berkes
Department of Statistics
Graz University of Technology
Steyrergasse 17
A-8010 Graz
Austria
and
A. Rényi Institute of Mathematics
Hungarian Academy of Sciences
P.O. Box 127
H-1364 Budapest
Hungary

L. Horváth
Department of Mathematics
University of Utah
155 South 1440 East
Salt Lake City, Utah 84112-0090
USA

P. Kokoszka
Department of Mathematics
and Statistics
Utah State University
3900 Old Main Hill
Logan, Utah 84322-3900
USA
E-mail: piotr@stat.usu.edu

Q.-M. Shao
Department of Mathematics
University of Oregon
Eugene, Oregon 97403-1222
USA
and
Department of Mathematics
Hong Kong University of Science
and Technology
Clear Water Bay, Kowloon
Hong Kong